\documentclass{amsart}
 
\RequirePackage{amsmath}
\RequirePackage{amssymb}
 
\theoremstyle{plain}
\newtheorem{theorem}{Theorem}[section]
\newtheorem{proposition}[theorem]{Proposition}
\newtheorem{lemma}[theorem]{Lemma}
\newtheorem{corollary}[theorem]{Corollary}

\newtheorem*{mainthm}{Theorem}

\theoremstyle{definition}

\newtheorem{remark}[theorem]{Remark}

\numberwithin{equation}{section}

\newcommand{\F}{\mathbf{F}}
\newcommand{\Fl}{\F_{\ell}}

\newcommand{\Kbar}{\bar{K}}
\newcommand{\Q}{\mathbf{Q}}
\newcommand{\Qbar}{\bar{\Q}}
\newcommand{\Qp}{\Q_{p}}
\newcommand{\Qpbar}{\Qbar_{p}}

\newcommand{\GQ}{G_{\Q}}
\newcommand{\GQS}{G_{\Q,S}}
\newcommand{\GQSl}{G_{\Q,S\cup\{\ell\}}}

\newcommand{\CC}{\mathcal{C}}
\newcommand{\lb}{\bar{\lambda}}
\newcommand{\m}{\mathfrak{m}}
\renewcommand{\O}{\mathcal{O}}
\newcommand{\p}{\mathfrak{p}}
\newcommand{\Z}{\mathbf{Z}}
\newcommand{\Zl}{\Z_{\ell}}

\newcommand{\Hf}{H^{1}_{f}}
\newcommand{\Hw}{H^{1}_{w}}

\newfont{\cyrr}{wncyr10}
\newcommand{\Sha}{\mbox{\cyrr Sh}}

\newcommand{\Klm}{K_{\lambda}}
\newcommand{\Klmbar}{\Kbar_{\lambda}}
\newcommand{\klm}{k_{\lambda}}
\newcommand{\klmbar}{\bar{k}_{\lambda}}
\newcommand{\Olm}{\O_{\lambda}}

\newcommand{\rhob}{\bar{\rho}}
\newcommand{\rhobflm}{\rhob_{f,\lambda}}

\newcommand{\rhoflm}{\rho_{f,\lambda}}

\newcommand{\chib}{\bar{\chi}}

\newcommand{\epl}{\varepsilon_{\ell}}
\newcommand{\eplb}{\bar{\varepsilon}_{\ell}}
\newcommand{\eplm}{\varepsilon_{\ell}}

\newcommand{\inj}{\hookrightarrow}
\newcommand{\surj}{\twoheadrightarrow}
\newcommand{\toi}{\overset{\simeq}{\longrightarrow}}
\newcommand{\too}{\longrightarrow}
\newcommand{\ts}{\textstyle}

\DeclareMathOperator{\Abs}{Red}
\DeclareMathOperator{\ad}{ad}
\DeclareMathOperator{\Cong}{Cong}
\DeclareMathOperator{\crys}{crys}
\DeclareMathOperator{\End}{End}
\DeclareMathOperator{\Frob}{Frob}
\DeclareMathOperator{\Gal}{Gal}
\DeclareMathOperator{\GL}{GL}
\DeclareMathOperator{\Hom}{Hom}
\DeclareMathOperator{\Obs}{Obs}
\DeclareMathOperator{\Sets}{Sets}

\DeclareMathOperator{\univ}{univ}

\newcommand{\adz}{\ad^{0}\!}

\begin{document}

\title[Explicit unobstructed primes]{Explicit unobstructed primes for modular deformation problems of
squarefree level}
\author{Tom Weston}
\email{weston@math.berkeley.edu}
\address{Department of Mathematics, University of California, Berkeley}
\thanks{Supported by an NSF postdoctoral fellowship}
\dedicatory{Dedicated to the memory of Arnold Ephraim Ross}

\begin{abstract}
Let $f$ be a newform of weight at least $2$ and squarefree level 
with Fourier coefficients
in a number field $K$.  We give explicit bounds, depending on congruences
of $f$ with other newforms, on the set of primes $\lambda$ of
$K$ for which the deformation problem associated to the mod $\lambda$
Galois representation of $f$ is obstructed.
We include some explicit examples.
\end{abstract}

\maketitle

\section{Introduction}

Let $f$
be a newform of weight $k \geq 2$, level $N$, and character $\omega$.
Let $K$ be the number field generated by the Fourier coefficients of $f$.
For any prime $\lambda$ of $K$ Deligne
has constructed a semisimple Galois representation
\[\rhobflm : \GQSl \to \GL_{2} \klm\]
over the residue field $\klm$ of $K$ at
$\lambda$; here $\GQSl$ is the Galois group
of the maximal extension of $\Q$ unramified outside the set
$S$ of places dividing $N\infty$ and the characteristic $\ell$ of $\klm$.
The representation $\rhobflm$ is absolutely irreducible for almost all
primes $\lambda$;
we write $\Abs(f)$ for the set of $\lambda$ such that
$\rhobflm$ is not absolutely irreducible.

Following Mazur,
we say that a prime $\lambda \notin \Abs(f)$ is an {\it obstructed prime}
for $f$ if the cohomology group
$H^{2}(\GQSl,\ad \rhobflm)$ of the adjoint representation of $\rhobflm$
is non-zero.  We write
$\Obs(f)$ for the set of such primes.  
The importance of this notion
rests on the fact that for $\lambda \notin
\Obs(f) \cup \Abs(f)$,
the universal deformation ring associated to
$\rhobflm$ is isomorphic to a power series ring in three variables over
the Witt vectors of $\klm$; see Section~\ref{s2} for details.

It was shown in
\cite{Weston2} that $\Obs(f)$ is finite for $f$ of weight $k \geq 3$.
In this paper
we obtain an explicit bound on $\Obs(f)$ in the case that the level $N$
of $f$ is squarefree.
We state our result here only for $N >1$;
see Section~\ref{s42} for the general statement (where we also allow
$k=2$ and $S$ non-minimal) and a partial converse.

\begin{mainthm}
Assume that $k \geq 3$ and that $N > 1$ is squarefree.  
Let $M$ denote the conductor of the Dirichlet character $\omega$.  Then
\[\Obs(f) 
\subseteq \bigl\{ \lambda \mid \ell \,;\, \ell \leq k+1 \text{~or~} \ell \mid
N\varphi(N){\ts \underset{p \mid \frac{N}{M}}{\prod}(p+1)} 
\bigr\} \cup \Cong(f)\]
with $\Cong(f)$ the set of congruence primes for $f$ (as defined
in Section~\ref{s41}) and $\varphi$ the Euler totient function.
\end{mainthm}

We note that the set
$\Cong(f)$ is computable using the results of \cite{Sturm} and
a tool such as \cite{Stein}.
It is not immediately clear to the author what form to expect the analogue
of this result for to take for $N$ not squarefree.

In Section~\ref{s2} we give a brief review of deformation theory and
use standard duality arguments to reduce the vanishing
of $H^{2}(\GQS,\ad \rhobflm)$ to the vanishing of certain local and
global cohomology groups.  The local groups are the subject of
Section~\ref{s3}; the computations rest on some simple cases of the
local Langlands correspondence.
In Section~\ref{s41} we use results of Hida (as refined in
\cite{Ghate}) to relate the global cohomology group to a certain Selmer
group studied by Diamond, Flach, and Guo.
The main results of the paper are proved in Section~\ref{s42}.
We give several explicit examples in Section~\ref{s5}.

It is a pleasure to thank Matthias Flach, Elena Mantovan, Robert Pollack, and
Ken Ribet for helpful conversations related to this paper.

\subsection*{Notation}

If $\rho : G \to \GL_{2}\!R$ 
is a representation of a group $G$ over a ring $R$,
we write $\ad \rho : G \to \GL_{4}\!R$ for the adjoint representation of $G$
on $\End(\rho)$ and $\adz \rho : G \to \GL_{3}\!R$ for the kernel of
the trace map from $\ad \rho$ to the trivial representation.
If $\rho : G \to \GL_{n}\!R$ is any representation, 
we write $H^{i}(G,\rho)$ for the cohomology group
$H^{i}(G,V_{\rho})$ with $V_{\rho}$ a free $R$-module of rank $n$ with
$G$-action via $\rho$.

We write $\GQ$ for the absolute Galois group of $\Q$.  We fix now and
forever embeddings $\Qbar \inj \Qpbar$ for each $p$, yielding injections
$G_{p} \inj \GQ$ with $G_{p}$ the absolute Galois group of $\Qp$.
We write $I_{p}$ for the inertia subgroup of $G_{p}$.
Let $\epl : \GQ \to \Zl^{\times}$ be the $\ell$-adic cyclotomic
character and let $\eplb : \GQ \to \Fl^{\times}$ be its reduction, the
mod $\ell$ Teichm\"uller character.
If $M$ is a $\Zl[\GQ]$-module, we write $M(1)$ for its first Tate
twist $M \otimes_{\Zl} \epl$.
If $S$ is a set of places of $\Q$ containing the infinite place, the
expression ``$p \in S$'' is to be interpreted as ``$p \in S - \{\infty\}$''.

\section{Obstructions} \label{s2}

\subsection{Deformation theory} \label{s21}

In this section we review the fundamentals of the deformation theory
of representations of profinite groups as in \cite{Mazur}.
Let $k$ be a finite field and
let $\CC$ denote the category of local rings which
are inverse limits of artinian local rings with residue field $k$;
a morphism $A \to B$ in $\CC$ is a continuous local homomorphism inducing the
identity map on residue fields.
Note that any ring $A$ in $\CC$ is canonically an algebra for the
Witt vectors $W(k)$ of $k$.

Let $G$ be a profinite group and fix an absolutely irreducible
continuous representation
\[\rhob : G \to \GL_{n}\!k\]
for some $n \geq 1$.  A {\it lifting} of $\rhob$ to a ring $A$ in $\CC$
is a continuous 
representation $\rho : G \to \GL_{n}\!A$ such that the composition
\[G \overset{\rho}{\too} \GL_{n}\!A \too \GL_{n}\!k\]
is equal to $\rho$.  Two liftings $\rho_{1},\rho_{2}$ of $\rhob$ are
said to be {\it strictly equivalent} if there is a matrix $M$ in the
kernel of $\GL_{n}\!A \to \GL_{n}\!k$
such that $\rho_{1} = M \cdot \rho_{2} \cdot M^{-1}$.  

A {\it deformation}
of $\rhob$ to $A$ is a strict equivalence class of liftings.  Let
\[D_{\rhob} : \CC \to \Sets\]
be the functor sending a ring $A$ to the set of deformations of $\rhob$ to
$A$.  The deformation functor $D_{\rhob}$ is representable by 
\cite[Section 1.2]{Mazur}; that is, there is a
ring $R_{\rhob}$ in $\CC$ (called the {\it universal deformation ring}
of $\rhob$) and an isomorphism of functors
\begin{equation} \label{eq:func}
D_{\rhob}(-) \cong \Hom_{\CC}(R_{\rhob},-).
\end{equation}
Note that via (\ref{eq:func}) the identity map on $R_{\rhob}$ corresponds
to a deformation
\[\rho^{\univ} : G \to \GL_{n}\!R_{\rhob}\]
of $\rhob$ to $R_{\rhob}$; this is the {\it universal deformation} of
$\rhob$, and the isomorphism (\ref{eq:func}) sends
$f : R_{\rhob} \to A$ to the deformation $f \circ \rho^{\univ}$ of
$\rhob$ to $A$.

The next proposition gives the fundamental connection between the deformation
problem $D_{\rhob}$ and the cohomology groups $H^{i}(G,\ad \rhob)$.
We say that $D_{\rhob}$ is {\it unobstructed}
if $H^{2}(G,\ad \rhob) = 0$.

\begin{proposition} \label{prop:defthy}
Assume that $H^{i}(G,\ad \rhob)$ is finite-dimensional over $k$ for
each $i$; set $d =dim_{k} H^{1}(G,\ad \rhob)$.
Then there exists a (non-canonical) surjection
\begin{equation} \label{eq:surj}
W(k)[[T_{1},\ldots,T_{d}]] \surj R_{\rhob}
\end{equation}
with kernel generated by at most $\dim_{k} H^{2}(G,\ad \rhob)$ elements.
In particular, if $D_{\rhob}$ is unobstructed, then
(\ref{eq:surj}) is an isomorphism.
\end{proposition}
\begin{proof}
This is proved in \cite[Section 1.6]{Mazur}.  The existence of the
surjection (\ref{eq:surj}) follows from an isomorphism
\[D_{\rhob}\bigl(k[\epsilon]/\epsilon^{2}\bigr) \toi H^{1}(G,\ad \rhob) \]
(sending a deformation $\rho$ to the cocycle $c_{\rho}$ such that
$\rho(g) = \rhob(g)(1+\epsilon \cdot c_{\rho}(g))$ for all $g \in G$) and
the interpretation of these groups as the tangent space of $R_{\rhob}$
via (\ref{eq:func}).  The statement about the kernel $J$ of (\ref{eq:surj})
follows from an injection
\[ \Hom(J,k) \inj H^{2}(G,\ad \rhob) \]
constructed using an obstruction two-cocycle measuring the failure of
$\rho^{\univ}$ to lift via (\ref{eq:surj}).
\end{proof}

The next lemma will be useful later in the paper.

\begin{lemma} \label{lemma:adjvan}
Let $\rhob : G \to \GL_{2}\!k$ be continuous and absolutely irreducible
and let $\chi : G \to k^{\times}$ be a character of
order at least $3$.  Then $H^{0}(G,\chi \otimes \ad \rhob)=0$.
\end{lemma}
\begin{proof}
If the image of $\rhob$ is dihedral, then the $G$-representation 
$\ad \rhob$ is the sum of the trivial character, a quadratic character,
and an irreducible two-dimensional representation of $G$.
If the image of $\rhob$ is not dihedral, then $\ad \rhob$ is the
sum of the trivial character and an irreducible three-dimensional
representation of $G$.  In either case the lemma follows since $\chi$ is
neither trivial nor quadratic.
\end{proof}

\subsection{Galois cohomology} \label{s22}

Let $k$ be a finite field of odd characteristic $\ell$.
We now apply the discussion of the previous section
to the case of a two-dimensional Galois representation over $k$.
Fix a finite set $S$ of places of $\Q$ including $\ell$ and the
infinite place.  Let $\Q_{S}$ denote the maximal extension of $\Q$ unramified
outside $S$; set $\GQS := \Gal(\Q_{S}/\Q)$.  Let
\[\rhob : \GQS \to \GL_{2}\!k\]
be continuous and absolutely irreducible.
We assume further that $\rhob$
is {\it odd} in the sense that the image of complex conjugation
has distinct eigenvalues.  In this section we study
the cohomology groups $H^{i}(\GQS,\ad \rhob)$.

\begin{lemma} \label{lemma:delta}
Each cohomology group $H^{i}(\GQS,\ad \rhob)$ is finite-dimensional over $k$
and 
\[\dim_{k} H^{1}(\GQS,\ad \rhob) - \dim_{k} H^{2}(\GQS,\ad \rhob) = 3.\]
\end{lemma}
\begin{proof}
The first statement is \cite[Corollary 4.15]{Milne}, while the second
is a straightforward calculation using Tate's global Euler
characteristic formula as in \cite[Section 1.10]{Mazur}.
\end{proof}

\begin{corollary} \label{cor:unob}
If $H^{2}(\GQS,\ad \rhob)=0$, then the universal deformation ring
$R_{\rhob}$ is (non-canonically) isomorphic to
$W(k)[[T_{1},T_{2},T_{3}]]$.
\end{corollary}

We will use global duality theorems of Poitou and Tate
to study $H^{2}(\GQS,\ad \rhob)$.
For a $k[\GQS]$-module $M$, define
\[\Sha^{1}(\GQS,M) := \ker\bigl(H^{1}(\GQS,M) \to \underset{p \in S}{\oplus}
H^{1}(G_{p},M)\bigr).\]

\begin{lemma} \label{lemma:d2}
One has
\[\dim_{k} H^{2}(\GQS,\ad \rhob) \leq
\dim_{k} \Sha^{1}(\GQS,\eplb \otimes \adz \rhob) + 
\underset{p \in S}{\ts \sum} \dim_{k} H^{0}(G_{p},\eplb \otimes \ad \rhob)\]
with equality if $\ell \neq 3$.
\end{lemma}
\begin{proof}
The trace pairing $\ad \rhob \otimes \ad \rhob \to k$ identifies
$\eplb \otimes \ad \rhob$ with the Cartier dual of $\ad \rhob$.  Thus
by \cite[Theorem 4.10]{Milne} there is an exact sequence
\begin{multline*}
0 \to H^{0}(\GQS,\eplb \otimes \ad \rhob) \to
\underset{p \in S}{\oplus} 
H^{0}(G_{p},\eplb \otimes \ad \rhob) \to \\
\Hom\bigl(H^{2}(\GQS,\ad \rhob),k\bigr) \to 
\Sha^{1}(\GQS,\eplb \otimes \ad \rhob) \to 0.
\end{multline*}
Since $\eplb \otimes \ad \rhob = \eplb \oplus (\eplb \otimes \adz \rhob)$
and $\Sha^{1}(\GQS,\eplb)$  vanishes by \cite[Lemma 10.6]{Weston},
the lemma follows from the exact sequence and Lemma~\ref{lemma:adjvan}.
\end{proof}

We will study the local terms $H^{0}(G_{p},\eplb \otimes \ad \rhob)$
in Section~\ref{s3}.  The global term
$\Sha^{1}(\GQS,\eplb \otimes \adz \rhob)$ is difficult to control
directly; instead we now relate it to a certain Selmer group, which in turn
is often computable using the results of Section~\ref{s41}.

Fix a totally ramified extension $K$ of the field of fractions of
$W(k)$.
The ring of integers $\O$ of $K$ lies in $\CC$; we write
$\m$ for its maximal ideal.
Let $\rho : \GQS \to \GL_{2}\!\O$ be a lifting of $\rhob$ to $\O$.
Let $V_{\rho}$ (resp.\ $A_{\rho}$) denote a three-dimensional $K$-vector space
(resp.\ $(K/\O)^{3}$) endowed with a $\GQS$-action via
$\adz \rho : \GQS \to \GL_{3}\!\O$.

Let $V$ (resp.\ $A$) denote either $V_{\rho}$ (resp.\ $A_{\rho}$) or else
its Tate twist.  For a prime $p$, define
\[\Hf(G_{p},V) := \begin{cases} H^{1}(G_{p}/I_{p},V^{I_{p}}) & p \neq \ell; \\
\ker\bigl(H^{1}(G_{p},V) \to H^{1}(G_{p},V \otimes B_{\crys})\bigr) & p = \ell;
\end{cases}\]
regarded as a $K$-subspace of $H^{1}(G_{p},V)$; here $B_{\crys}$ is
the crystalline period ring of Fontaine.  Let
$\Hf(G_{p},A)$ denote the image of $\Hf(G_{p},V)$ under the
pushforward from $H^{1}(G_{p},V)$ to $H^{1}(G_{p},A)$.
For $M$ denoting either of $V$ or $A$, the {\it Selmer group} of $M$ is
defined by
\[\Hf(\GQ,M) := \bigl\{ c \in H^{1}(\GQ,M) \,;\, c|_{G_{p}} \in \Hf(G_{p},M)
\text{~for all~}p \bigr\}.\]

Following \cite[Section 7]{DFG},
we will also need a slight variant of this construction.  Define
\begin{gather*}
\Hw(G_{p},A) := \begin{cases} H^{1}(G_{p}/I_{p},A^{I_{p}}) & p \neq \ell; \\
\Hf(G_{\ell},A) & p = \ell; \end{cases} \\
H^{1}_{\emptyset}(\GQ,A) := \bigl\{ c \in H^{1}(\GQ,A) \,;\, c|_{G_{p}} \in 
\Hw(G_{p},A) \text{~for all~}p \bigr\}.
\end{gather*}
Clearly one has $\Hf(G_{p},A) \subseteq \Hw(G_{p},A)$ for all $p$, so that
\begin{equation} \label{eq:comp}
\Hf(\GQ,A) \subseteq H^{1}_{\emptyset}(\GQ,A).
\end{equation}
In fact, this inclusion is an equality if $A^{I_{p}}$ is divisible for
all $p \neq \ell$.

\begin{lemma} \label{lemma:ineq}
Assume that $\ell > 3$ and $\Hf(\GQ,V_{\rho}) = \Hf(\GQ,V_{\rho}(1))=0$.  Then
\[\dim_{k} \Sha^{1}(\GQS,\eplb \otimes \adz \rhob) \leq
\dim_{k} H^{1}_{\emptyset}(\GQS,A_{\rho})[\m].\]
\end{lemma}
\begin{proof}
Since $\rho$ is a lifting of $\rhob$, the $k[\GQS]$-module
$A_{\rho}(1)[\m]$ is a realization of $\eplb \otimes \adz \rhob$.
We thus obtain a natural map
\begin{equation} \label{eq:inj}
H^{1}(\GQ,\eplb \otimes \adz\rhob) = 
H^{1}\bigl(\GQ,A_{\rho}(1)[\m]\bigr) \to H^{1}\bigl(\GQ,A_{\rho}(1)\bigr)
\end{equation}
which is injective by Lemma~\ref{lemma:adjvan}.
The image of $\Sha^{1}(\GQS,\eplb \otimes \adz \rhob)$ under (\ref{eq:inj})
is easily seen to
lie in $\Hf(\GQ,A_{\rho}(1))$, so that we obtain an injection
\begin{equation} \label{eq:inj2}
\Sha^{1}(\GQS,\eplb \otimes \adz \rhob) \inj
\Hf\bigl(\GQ,A_{\rho}(1)\bigr).
\end{equation}
By \cite[Theorem 1]{Flach1}, the latter group is (non-canonically)
isomorphic to $\Hf(\GQ,A_{\rho})$ (see \cite[Proposition 2.2]{Weston2};
this also uses the assumption on the vanishing of the Selmer group of
$V_{\rho}$ and $V_{\rho}(1)$).
The lemma thus follows from (\ref{eq:inj2}) and (\ref{eq:comp}).  
\end{proof}

\begin{remark}
The only difficulty in analyzing the failure of (\ref{eq:inj2}) to
be an isomorphism on $\m$-torsion is the determination of the image of
the restriction map
$$\Hf\bigl(\GQ,A_{\rho}(1)\bigr)[\m] \to \Hf\bigl(G_{\ell},A_{\rho}(1)
\bigr)[\m].$$
Unfortunately, this question appears to be quite difficult in general.
\end{remark}

\section{Local invariants} \label{s3}

Let $f = \sum a_{n}q^{n}$ be a newform of weight $k \geq 2$, squarefree
level $N$, and character $\omega$.  Let $K$ be the number field generated
by the Fourier coefficients
$a_{n}$ of $f$.  For any prime $\lambda$ of $K$, Deligne has constructed
a continuous $\lambda$-adic Galois representation
\[\rhoflm : \GQ \to \GL_{2}\! \Klm.\]
This representation is unramified at $p \nmid N\ell$ (with $\ell$ the
characteristic of the residue field $\klm$ of $\Klm$) and for such
$p$ the trace (resp.\ the determinant)
of the image of an arithmetic Frobenius element $\Frob_{p}$
under $\rhoflm$ is equal to $a_{p}$ (resp.\ $p^{k-1}\omega(p)$).

As usual we identify $\omega : (\Z/N\Z)^{\times} \to \mu_{\varphi(N)}$ 
with a Galois
character via the canonical isomorphism $\Gal(\Q(\mu_{N})/\Q) \cong
(\Z/N\Z)^{\times}$; the determinant of $\rhoflm$ is then
$\eplm^{k-1}\omega$.
Let $M$ denote the conductor of $\omega$ and let
$\omega_{0} : (\Z/M\Z)^{\times} \to  \mu_{\varphi(M)}$
be the associated primitive Dirichlet character.
Then $\omega$ is ramified at $p$
if and only if $p$ divides $M$, in which case the restriction of
$\omega$ to the inertia group 
$I_{p}$ is a non-trivial character taking values in $\mu_{p-1}$.

For the remainder of this section we fix a prime $\lambda$ of $K$
dividing a rational prime $\ell$.
Let
\[\rhobflm : \GQ \to \GL_{2}\!\klm\]
be the semisimple reduction of $\rhoflm$; this is well-defined independent
of any choice of integral model of $\rhoflm$.  We are interested in the
local invariants $H^{0}(G_{p},\eplb \otimes \ad \rhobflm)$ for all
primes $p$.  As
\[\eplb \otimes \ad \rhobflm \cong \eplb \oplus (\eplb \otimes \adz
\rhobflm)\]
and
\[H^{0}(G_{p},\eplb) \neq 0 \,\, \Leftrightarrow \,\, 
p \equiv 1 \pmod{\ell},\]
we will restrict our attention below to the case that $\ell$ does not
divide $p-1$ and to the study of
$H^{0}(G_{p},\eplb \otimes \adz \rhobflm)$.

In the analysis below
we make use of the local Langlands correspondence and the
compatibility results completed in \cite{Carayol}.
Rather than review these results
in detail, we will only recall the consequences we need; see
\cite{Weston2} for more details and references.

\subsection{$\boldsymbol{p \nmid N\ell}$} \label{s31}

Fix $\alpha_{p},\beta_{p} \in \Kbar$ with $\alpha_{p}+\beta_{p} = a_{p}$
and $\alpha_{p}\beta_{p} = p^{k-1}\omega(p)$.
In this case, we have
\[\rhoflm|_{G_{p}} \otimes \Klmbar \cong \chi_{1} \oplus \chi_{2}\]
where the $\chi_{i} : G_{p} \to \Klmbar^{\times}$
are unramified characters with
\begin{equation} \label{eq:char}
\chi_{1}(\Frob_{p}) = \alpha_{p}; \qquad
\chi_{2}(\Frob_{p}) = \beta_{p}.
\end{equation}
We write $\chib_{i} : G_{p} \to \klmbar^{\times}$ for the reduction of
$\chi_{i}$.

\begin{lemma} \label{lemma:unramps}
Assume $p \nmid N\ell$ and $p \not\equiv 1 \pmod{\ell}$.  Then
$H^{0}(G_{p},\eplb \otimes \ad \rhobflm) \neq 0$
if and only if
$a_{p}^{2} \equiv (p+1)^{2}p^{k-2}\omega(p) \pmod{\lambda}$.
\end{lemma}
\begin{proof}
Since the existence of eigenvectors with $\klm$-rational eigenvalues is
invariant under base extension, the existence of $G_{p}$-invariants
in $\eplb \otimes \adz \rhobflm$ is equivalent to the existence of
$G_{p}$-invariants in
\[\bigl(\eplb \otimes \adz \rhobflm|_{G_{p}} \bigr)
\otimes \klmbar \cong \eplb \oplus 
\eplb\chib_{1}^{\vphantom{-1}}\chib_{2}^{-1} \oplus
\eplb\chib_{1}^{-1}\chib_{2}^{\vphantom{-1}}.\]
As $p \not\equiv 1 \pmod{\ell}$,  this has
non-trivial $G_{p}$-invariants if and only if one of the characters
$\eplb\chib_{1}^{\vphantom{-1}}\chib_{2}^{-1}$,
$\eplb\chib_{1}^{-1}\chib_{2}^{\vphantom{-1}}$
is trivial.  By (\ref{eq:char}) this occurs if and only if
\begin{alignat*}{2}
\frac{\alpha_{p}}{\beta_{p}} &\equiv p^{\pm 1} & &\pmod{\lambda}. \\
\intertext{This in turn is equivalent to}
\frac{\alpha_{p}}{\beta_{p}} + \frac{\beta_{p}}{\alpha_{p}} &\equiv
p + \frac{1}{p} & &\pmod{\lambda}\\
\frac{(\alpha_{p}+\beta_{p})^{2}}{\alpha_{p}\beta_{p}} &\equiv
\frac{(p+1)^{2}}{p} & &\pmod{\lambda}\\
a_{p}^{2} &\equiv (p+1)^{2}p^{k-2}\omega(p) & &\pmod{\lambda}
\end{alignat*}
as claimed.
\end{proof}

\subsection{$\boldsymbol{p \mid M}$, $\boldsymbol{p \neq \ell}$} \label{s32}

In this case the $p$-component $\pi_{p}$ of
the automorphic representation associated to $f$ has conductor $1$ and
ramified central character.  It follows that $\pi_{p}$ is
a principal series representation associated to one ramified
character and one unramified character.
On the Galois side, this translates to
\[\rhoflm|_{G_{p}} \otimes \Klmbar \cong \chi_{1} \oplus \chi_{2}\]
for continuous characters $\chi_{i} : G_{p} \to \Klmbar^{\times}$ with
$\chi_{1}$ ramified and $\chi_{2}$ unramified.
Since $\rhoflm|_{G_{p}}$ has determinant $\eplm^{k-1}\omega|_{G_{p}}$,
we have $\chi_{1}\chi_{2} = \eplm^{k-1}\omega|_{G_{p}}$.  In particular
$\chi_{1}|_{I_{p}} = \omega|_{I_{p}}$ is a non-trivial character taking
values in $\mu_{p-1}$.
If $p \not \equiv 1 \pmod{\ell}$, then $\mu_{p-1}$ injects
into $\klm^{\times}$ and consequently the reduction $\bar{\chi}_{1} : G_{p} \to
\klmbar^{\times}$ is still ramified at $p$.

\begin{lemma} \label{lemma:ramps}
Assume $p \mid M$, $p \neq \ell$, and $p \not\equiv 1 \pmod{\ell}$.  Then
$H^{0}(G_{p},\eplb \otimes \ad \rhobflm) = 0$.
\end{lemma}
\begin{proof}
As in Lemma~\ref{lemma:unramps}, it suffices to show that
the two characters
$\eplb\chib_{1}^{\vphantom{-1}}\chib_{2}^{-1}$,
$\eplb\chib_{1}^{-1}\chib_{2}^{\vphantom{-1}}$
are non-trivial.  Since 
$\eplb$ and $\chib_{2}$ are unramified at $p$ while
$\chib_{1}$ is ramified at $p$, this is clear.
\end{proof}

\subsection{$\boldsymbol{p \mid \frac{N}{M}}$, $\boldsymbol{p \neq \ell}$} 
\label{s33}

In this case $\pi_{p}$ has conductor $1$ and unramified central
character.  It follows that $\pi_{p}$ is the special representation
associated to an unramified character.  This means that there
exists an unramified character $\chi : G_{p} \to \Klmbar^{\times}$
such that
\begin{equation} \label{eq:spec}
\rhoflm|_{G_{p}} \otimes \Klmbar \cong 
\left(\begin{array}{cc} \eplm\chi & * \\ 0 & \chi \end{array}\right)
\end{equation}
with the upper right corner ramified.

\begin{lemma} \label{lemma:specss}
Assume $p \mid \frac{N}{M}$, $p \neq \ell$, and
$p^{2} \not\equiv 1 \pmod{\ell}$.  Then
$H^{0}(G_{p},\eplb \otimes \ad \rhobflm) \neq 0$
if and only if $\rhobflm$ is unramified at $p$.
\end{lemma}
\begin{proof}
Since $p^{2} \not\equiv 1 \pmod{\ell}$, 
by \cite[Lemma 5.1]{Weston2} we have
\[\rhobflm|_{G_{p}} \otimes \klmbar 
\cong \left(\begin{array}{cc} \eplb\chib & \nu
\\ 0 & \chib \end{array}\right)\]
for some $\nu : G_{p} \to \klmbar$; in fact, one checks directly
that $\chib^{-1}\nu$
is naturally an element of $H^{1}(G_{p},\klmbar(1))$.  Since $\eplb$ and
$\chib$ are unramified, $\rhobflm|_{G_{p}}$ is unramified if and only if
$\chib^{-1}\nu$ is unramified.  
However, since $p \not\equiv 1 \pmod{\ell}$
every non-zero element of $H^{1}(G_{p},\klmbar(1))$ is ramified.
We conclude that 
$\rhobflm|_{G_{p}}$ is unramified if and only if it is
semisimple.  \cite[Lemma 5.2]{Weston2} now completes the proof.
\end{proof}

\begin{lemma} \label{lemma:spec}
Assume $p \mid \frac{N}{M}$, $p \neq \ell$,
$p^{2} \not\equiv 1 \pmod{\ell}$, and $\rhobflm$
absolutely irreducible.  Then
$H^{0}(G_{p},\eplb \otimes \ad \rhobflm) \neq 0$ if and only if
there exists a newform $f'$, of weight $k$ and level dividing $\frac{N}{p}$,
such that $\rhob_{f,\lb} \cong \rhob_{f',\lb}$
for some prime $\lb$ of $\Qbar$ over $\lambda$.
\end{lemma}

Here by $\rhob_{f,\lb}$ (resp.\ $\rhob_{f',\lb}$)
we mean $\rhob_{f,\lambda} \otimes \klmbar$ (resp.\ 
$\rhob_{f',\lambda'} \otimes \bar{k}'_{\lambda'}$ with $\lambda'$ the
intersection of $\lb$ with the field $K'$ of Fourier coefficients
of $f'$ and with $k'_{\lambda'}$ the residue field of $K'$ at $\lambda'$.)

\begin{proof}
By \cite[(B) of p.\ 221]{Edixhoven}, the existence of such an $f'$ is
equivalent to $\rhobflm$ being unramified at $p$.  Thus the lemma follows
from Lemma~\ref{lemma:specss}.
\end{proof}

\begin{remark} \label{rmk:spec}
If one further assumes that $p' \not\equiv 1 \pmod{\ell}$ for all $p'$ 
dividing $N$, then the newform $f'$ of Lemma~\ref{lemma:spec}
must have level a multiple of $M$ and
character lifting $\omega_{0}$, so that $\lambda$ is a congruence prime
for $f$ of level dividing $\frac{N}{p}$ in the terminology
of Section~\ref{s41}.  Indeed, $\rhob_{f',\lambda}$ is
isomorphic to $\rhob_{f,\lambda}$ and thus has determinant
$\eplb^{k-1}\bar{\omega}$; therefore the character
$\omega'$ of $f'$
must have reduction equal to $\bar{\omega}$.  However, since
$p \not\equiv 1 \pmod{\ell}$ for all $p$ dividing $N$, the only such
characters of conductor dividing $N$
are those which lift $\omega_{0}$.  Thus $f'$ must have
level divisible by $M$ and character lifting $\omega_{0}$, as claimed.
\end{remark}

\subsection{$\boldsymbol{p = \ell}$, $\boldsymbol{\ell \nmid N}$} \label{s34}

We now give some mild improvements 
on the results of \cite[Section 4]{Weston2} on the vanishing
of $H^{0}(G_{\ell},\eplb \otimes \ad \rhobflm)$.
Recall that $f = \sum a_{n}q^{n}$ 
is said to be {\it ordinary} (resp.\ {\it supersingular})
at $\lambda$ if $v_{\lambda}(a_{\ell}) = 0$
(resp.\ $v_{\lambda}(a_{\ell}) > 0$), with $v_{\lambda}$ the
$\lambda$-adic valuation.  If $f$ is ordinary at $\lambda$,
then the semisimplification of $\rhoflm|_{I_{\ell}} \otimes \Klmbar$ is
isomorphic to $\eplm^{k-1} \oplus 1$,
while if $f$ is supersingular at $\lambda$, then $\rhobflm|_{G_{\ell}}$
is absolutely irreducible.  (This all follows from the discussion of
\cite[pp.\ 214--215]{Edixhoven}, for example.)

\begin{lemma} \label{lemma:ord}
Assume $\ell \nmid N$.
If $f$ is ordinary at $\lambda$ and 
$H^{0}(G_{\ell},\eplb \otimes \ad \rhobflm) \neq 0$, then
$k \equiv 0,2 \pmod{\ell - 1}$.
\end{lemma}
\begin{proof}
It suffices to prove the corresponding result for the 
$I_{\ell}$-invariants of the semisimplification
of $(\eplb \otimes \ad \rhobflm) \otimes \klmbar$.  
By the above discussion this semisimplification is isomorphic to
\[\eplb \oplus \eplb \oplus
\eplb^{k} \oplus \eplb^{2-k}.\]
Since $\eplb$ has order $\ell-1$, the lemma follows.
\end{proof}

Note that the above lemma is vacuous in the case of weight $2$.

\begin{lemma} \label{lemma:ss}
Assume $\ell \nmid N$.
If $f$ is supersingular at $\lambda$ and $\ell > 3$,
then $H^{0}(G_{\ell},\eplb \otimes \ad \rhobflm)=0$.
\end{lemma}
\begin{proof}
As $\rhobflm|_{G_{\ell}}$ is absolutely irreducible, this is immediate
from Lemma~\ref{lemma:adjvan}.
\end{proof}

\section{Global results} \label{s4}

We continue with a newform $f=\sum a_{n}q^{n}$ 
of weight $k$, squarefree level $N$, and
character $\omega$ of conductor $M$ 
as in Section~\ref{s3}.
Let $\O$ be the ring of integers of the field $K$ of Fourier coefficients
of $f$.  For each prime $\lambda$ of $K$, let $V_{\rho,\lambda}$ be
a three-dimensional $\Klm$-vector space with $\GQ$-action via
$\adz \rhoflm$.  Fix a $\GQ$-stable $\Olm$-lattice $T_{\rho,\lambda}$ in
$V_{\rho,\lambda}$ and set $A_{\rho,\lambda} := V_{\rho,\lambda}/
T_{\rho,\lambda}$.  
In general the $\klm[\GQ]$-module $A_{\rho,\lambda}[\lambda]$ 
need not agree with the semisimple
reduction $\rhobflm$; however, these two representations must be
isomorphic when
$\rhobflm$ is absolutely irreducible, which is the only case we will
consider below.

\subsection{Congruences and Selmer groups} \label{s41}

The purpose of this section is to explain how the results of
\cite{Hida} (as refined in \cite{Ghate})
and \cite{DFG} relate adjoint Selmer groups with congruences
of modular forms.  Let $d$ be a divisor of $N$ which is divisible by $M$.
We say that a 
prime $\lambda$ of $K$ is a {\it congruence prime of level $d$} for $f$ if
there exists a newform $f'$ of weight $k$ and level $d$ such that:
\begin{enumerate}
\item $f'$ has character lifting $\omega_{0}$;
\item $f'$ is not a Galois conjugate of $f$;
\item $\rhob_{f,\lb} \cong \rhob_{f',\lb}$
for some prime $\lb$ of $\Qbar$ above $\lambda$.
\end{enumerate}
(Of course, by the Cebatorev density theorem
the last condition is equivalent to a congruence
$a_{p}(f) \equiv a_{p}(f') \pmod{\lb}$ for all primes $p$
not dividing $N$.)

We say that a congruence prime $\lambda$ of level $d$ for $f$ is 
{\it proper} (resp.\ 
{\it strict}) if $d<N$ (resp.\ $d=N$).
Let $\Cong(f)$ (resp.\ $\Cong_{<N}(f)$, resp.\ $\Cong_{N}(f)$) denote
the set of congruence primes (resp.\ proper congruence primes,
resp.\ strict congruence primes) for $f$.

We need a simple lemma.

\begin{lemma} \label{lemma:dfg}
Let $\lambda$ be a prime of $K$ dividing a rational prime $\ell$ not
dividing $N$.
Assume that $N > 1$ and that $\rhobflm$ is ramified at some
$p$ dividing $N$.  If $\rhobflm$ is absolutely irreducible, then
$\rhobflm|_{G_{F}}$ is absolutely irreducible as well, where
$F = \Q(\sqrt{(-1)^{(\ell-1)/2}\ell})$.
\end{lemma}
\begin{proof}
As in \cite[Lemma 7.14]{DFG}, if $\rhobflm|_{G_{F}}$ is absolutely reducible,
then $\rhobflm$ is induced from a character of $G_{F}$.
In particular, it follows 
that the conductor $N'$ of $\rhobflm$ (in the sense of \cite{Edixhoven})
is a square.  However, $N'$ must also divide the level $N$ of $f$;
since $N'$ is non-trivial by hypothesis and
$N$ is squarefree, this is impossible.
\end{proof}

\begin{proposition} \label{prop:hida}
Let $\lambda$ be a prime of $K$ dividing the rational prime $\ell$.
Assume that:
\begin{enumerate}
\item \label{a1} $\rhobflm$ is absolutely irreducible;
\item \label{a2} $\ell > k$;
\item \label{a3} Either $N > 1$ or $\ell \nmid (2k-3)(2k-1)$;
\item \label{a4} $\ell \nmid N$;
\item \label{a5} $\ell \nmid \varphi(N)$ (that is,
$p \not\equiv 1 \pmod{\ell}$ for all $p \mid N$);
\item \label{a6} $\rhobflm$ is ramified at $p$ for all $p \mid 
\frac{N}{M}$;
\end{enumerate}
Then $H^{1}_{\emptyset}(\GQ,A_{\rho,\lambda}) \neq 0$ if and only if
$\lambda \in \Cong_{N}(f)$.
\end{proposition}
\begin{proof}
Conditions (\ref{a4})--(\ref{a6}) guarantee that
$A_{\rho,\lambda}$ is minimally ramified in the sense of
\cite[Section 3]{Diamond2}.  Using (\ref{a1}), (\ref{a2}), (\ref{a4}),
and Lemma~\ref{lemma:dfg} (or (\ref{a3}) and 
\cite[Lemma 7.14]{DFG} for $N=1$), 
we may apply \cite[Theorem 7.15]{DFG} to conclude that
\begin{equation} \label{eq:len}
\text{length}_{\Olm} 
H^{1}_{\emptyset}(\GQ,A_{\rho,\lambda}) = v_{\lambda}(\eta_{f}^{\emptyset});
\end{equation}
here $\eta_{f}^{\emptyset}$ is the fractional ideal of $K$ defined
in \cite[Section 6.4]{DFG} and $v_{\lambda}$ is the $\lambda$-adic valuation.

By definition,
the ideal $\eta_{f}^{\emptyset}$ is generated
by the discriminant $d(L_{f}(\O_{K}))$ of \cite[proof of Theorem 5]{Ghate},
which in turn equals the square of the
algebraic special value of the adjoint $L$-function of $f$:
\begin{equation} \label{eq:disc}
 d\bigl(L_{f}(\O_{K})\bigr) = \left(\frac{W(f)\Gamma(1,\ad f)L(1,\ad f)}
{\Omega(f,+)\Omega(f,-)}\right)^{2}.
\end{equation}
(All of this is only true
up to factors of primes violating (\ref{a1})--(\ref{a4}).)
In particular, (\ref{eq:len}) implies
that $H^{1}_{\emptyset}(\GQ,A_{\rho,\lambda})$
is non-zero if and only if
(\ref{eq:disc}) has positive $\lambda$-adic valuation.
By \cite[Theorems 1 and 2]{Ghate},
the latter condition is equivalent to the existence of a newform $f'$
of weight $k$ and level dividing $N$, not Galois conjugate to $f$, such
that $\rhob_{f,\lb} \cong \rhob_{f',\lb}$ for some prime $\lb$ above
$\lambda$.

It remains to show that $f'$ has level $N$ and character $\omega$.
Since $\rhob_{f',\lb}$ has determinant
$\eplb^{k-1}\bar{\omega}$ and $\mu_{\varphi(N)}$ injects into
$\klm^{\times}$ (by (\ref{a5})), $f'$ has level divisible by $M$
and character lifting
$\omega_{0}$.  Hypothesis (\ref{a6}) guarantees that
$\rhob_{f',\lb}$ is ramified at all $p \mid \frac{N}{M}$ as well,
so that $f'$ must in fact have level $N$.
\end{proof}

\subsection{Vanishing of cohomology} \label{s42}

Let $S$ be a finite set of places of $\Q$ containing all places 
dividing $N\infty$; let $N_{S}$ denote the product of all primes in $S$.
Fix a prime $\lambda$ of $K$ dividing a rational prime $\ell$.
We are now in a position to compute $H^{2}(\GQSl,\ad \rhobflm)$.

\begin{theorem} \label{thm:van}
Assume that $\rhobflm$ is absolutely irreducible and $\ell > 3$.  If
\begin{equation} \label{eq:h2}
H^{2}(\GQSl,\ad \rhobflm) \neq 0,
\end{equation}
then one of the following holds:
\begin{enumerate}
\item \label{b1} $\ell \leq k$;
\item \label{b2} $\ell \mid N$;
\item \label{b3} $\ell \mid \varphi(N_{S})$;
\item \label{b3.5} $\ell \mid p+1$ for some $p \mid \frac{N}{M}$;
\item \label{b4} $a_{p}^{2} \equiv (p+1)^{2}p^{k-2}\omega(p) \pmod{\lambda}$
for some $p \mid \frac{N_{S}}{N}$, $p \neq \ell$;
\item \label{b5} $\ell = k + 1$ and $f$ is ordinary at $\lambda$;
\item \label{b6} $k=2$ and $a_{\ell}^{2} \equiv \omega(\ell) \pmod{\lambda}$;
\item \label{b7} $N=1$ and $\ell \mid (2k-3)(2k-1)$;
\item \label{b8} $\lambda \in \Cong(f)$.
\end{enumerate}
\end{theorem}

Using Lemma~\ref{lemma:d2} and
the results of Sections~\ref{s3} and~\ref{s41}, the reader should
have little difficulty in detecting the source of each of the conditions
above.  We shall nevertheless endeavor to give a complete proof.

\begin{proof}
If (\ref{eq:h2}) holds, then
Lemma~\ref{lemma:d2} implies that either 
\begin{equation} \label{eq:pos1}
H^{0}(G_{p},\eplb \otimes \ad \rhobflm) \neq 0
\end{equation}
for some $p \in S \cup \{\ell\}$ or
\begin{equation} \label{eq:pos2}
\Sha^{1}(\GQS,\eplb \otimes \ad \rhobflm) \neq 0.
\end{equation}

Suppose first that (\ref{eq:pos1}) holds for a prime $p \in S \cup \{\ell\}$;
we may assume $\ell \nmid N$ by (\ref{b2}).
If $H^{0}(G_{p},\eplb) \neq 0$, then $\ell$ divides $p-1$ which in turn
divides $\varphi(N_{S})$, so that
(\ref{b3}) holds.  We may thus assume that $p \not\equiv 1 \pmod{\ell}$ and
\[H^{0}(G_{p},\eplb \otimes \adz \rhobflm) \neq 0.\]
By Lemma~\ref{lemma:ramps} we know that $p$ does not divide $M$.
If $p$ divides $\frac{N}{M}$, then by Lemma~\ref{lemma:spec} and
Remark~\ref{rmk:spec} one of (\ref{b3.5}) or (\ref{b8}) holds, while if
$p$ does not divide $N\ell$, then
Lemma~\ref{lemma:unramps} implies that (\ref{b4}) must hold.
Finally, if $p = \ell$ and $k > 2$,
then Lemmas~\ref{lemma:ord} and~\ref{lemma:ss} force
(\ref{b1}) or (\ref{b5}) to hold; if $k = 2$,
then \cite[Proposition 4.4]{Weston2} forces
(\ref{b6}) to hold.

It remains to consider the case that (\ref{eq:pos2}) holds, 
(\ref{eq:pos1}) does not hold for 
any $p \in S \cup \{\ell\}$, and none of (\ref{b1})--(\ref{b7}) hold.
Then by \cite[Theorem 8.2]{DFG}
\[ \Hf(\GQ,V_{f,\lambda}) = \Hf(\GQ,V_{f,\lambda}(1)) = 0.\]
Lemma~\ref{lemma:ineq} and (\ref{eq:pos2}) thus imply that
\[ H^{1}_{\emptyset}(\GQ,A_{\rho,\lambda}) \neq 0. \]
Since $H^{0}(G_{p},\eplb \otimes \ad \rhobflm) = 0$ for all $p$ dividing
$\frac{N}{M}$, Lemma~\ref{lemma:specss} implies that $\rhobflm$ is ramified
at all such $p$.  Proposition~\ref{prop:hida} now
applies to show that $\lambda \in \Cong_{N}(f)$.
Thus (\ref{b8}) holds, completing the proof.
\end{proof}

\begin{corollary} \label{cor:def}
If $\rhobflm$ is absolutely irreducible, $\ell > 3$, and $\lambda$
does not satisfy (\ref{b1})--(\ref{b8}), then
\[ R_{\rhobflm} \cong W(\klm)[[T_{1},T_{2},T_{3}]]. \]
\end{corollary}
\begin{proof}
This follows immediately from Theorem~\ref{thm:van} and
Corollary~\ref{cor:unob}.
\end{proof}

We also obtain the following partial converse to Theorem~\ref{thm:van}.

\begin{theorem} \label{thm:vanc}
Assume that $\rhobflm$ is absolutely irreducible.
Suppose that $\ell > 3$ and one of the following holds:
\begin{enumerate}
\item \label{m1} $\ell \mid \varphi(N_{S})$;
\item \label{m2} $a_{p}^{2} \equiv (p+1)^{2}p^{k-2}\omega(p) \pmod{\lambda}$
for some $p \mid \frac{N_{S}}{N}$, $p \neq \ell$;
\item \label{m3} $\lambda$ is a congruence prime for $f$ of 
level dividing $\frac{N}{p}$ for some
$p \mid \frac{N}{M}$, $\ell \nmid p(p+1)$;
\item \label{m4} $k=2$, $\ell \nmid N$, 
and $a_{\ell}^{2} \equiv \omega(\ell) \pmod{\lambda}$.
\end{enumerate}
Then $H^{2}(\GQSl,\eplb \otimes \ad \rhobflm) \neq 0$.
\end{theorem}
\begin{proof}
By Lemma~\ref{lemma:d2} it suffices to show that these conditions
guarantee that $H^{0}(G_{p},\eplb \otimes \ad \rhobflm) \neq 0$ for
some $p \in S$.  If (\ref{m1}) holds, then $H^{0}(G_{p},\eplb) \neq 0$
for some $p \in S$, so that this is clear.  If (\ref{m2}) holds, then
by Lemma~\ref{lemma:unramps} we have $H^{0}(G_{p},\eplb \otimes \adz 
\rhobflm) \neq 0$.  If (\ref{m1}) does not hold but
(\ref{m3}) does hold, then Lemma~\ref{lemma:spec} implies that
$H^{0}(G_{p},\eplb \otimes \adz\rhobflm) \neq 0$.
Finally, if (\ref{m4}) holds, then the proof of \cite[Proposition 4.4]{Weston2}
shows that $H^{0}(G_{\ell},\eplb \otimes \ad \rhobflm) \neq 0$.
\end{proof}

\section{Examples} \label{s5}

In this section we use the data of \cite{Stein} to bound the
obstructed primes for the
deformation problems associated to a few specific modular forms.
Of course, the most interesting aspect of these computations are
the determination of congruences between newforms.
Using \cite{Stein} we can check such congruences on the
$p^{\text{th}}$ Fourier coefficients for all $p < 1000$;
by the results of \cite{Sturm} these checks are more than sufficient
to prove that these congruences actually exist in our examples.
We will not comment further on this issue.

For a modular form $f$, we let
$\Abs(f)$ denote the set of primes $\lambda$ of $K$ such that
$\rhobflm$ is absolutely reducible.
We recall the following well-known facts regarding $\Abs(f)$; 
see \cite[Lemma 7.13]{DFG} for example.

\begin{lemma} \label{lemma:red}
Let $f = \sum a_{n}q^{n}$ 
be a newform of weight $k$ and level $N$ with coefficient field $K$.
Let $\lambda$ be a prime of $K$ dividing a rational prime $\ell$.
Suppose that $\lambda \in \Abs(f)$, so that
\[\rhobflm \otimes \klmbar \cong \chi_{1} \oplus \chi_{2}\]
for characters $\chi_{1},\chi_{2} : \GQ \to
\klmbar^{\times}$.  If $\ell$ does not divide $N$, then each $\chi_{i}$ has
conductor dividing $N\ell$.  If also $\ell > k$, then
one of the $\chi_{i}$ has conductor dividing $N$, so that
\[a_{p} \equiv p^{k-1} + 1 \pmod{\lambda}\]
for all $p \equiv 1 \pmod{N}$.
\end{lemma}

In practice one uses the second condition to bound the set
$\Abs(f)$ and the first condition to check each remaining $\lambda$
not dividing $N$.
For a prime $\lambda$ dividing $N$, one can still check that
$\rhobflm$ is absolutely reducible, but it is much more difficult to show
that $\rhobflm$ is absolutely irreducible; we will make no attempt to deal
with this case below.

For a finite set of places $S$ containing all places dividing $N\infty$,
we let $\Obs_{S}(f)$ denote the set of $\lambda \notin \Abs(f)$ such that
$$H^{2}(\GQSl,\epl \otimes \ad \rhobflm) \neq 0,$$
or equivalently such that
the deformation problem associated to 
\[\rhobflm : \GQSl \to \GL_{2}\! \klm\] 
is obstructed.
We simply write $\Obs(f)$ for
$\Obs_{\{p\mid N\infty\}}(f)$.

In the interests of space, we make the following notational conventions.
Fix a quadratic extension $K$ of $\Q$ and let $p$ be a rational prime.
If $p$ ramifies in $K$, then we simply write $\p_{p}$ for the prime of
$K$ above $p$.  If $p$ splits, then we will write $\p_{p}$ and $\bar{\p}_{p}$
for the two primes of $K$ above $p$, at least
when it is not important to distinguish between them.





\subsection{Weight $\boldsymbol{12}$, level $\boldsymbol{5}$, trivial
character}

There are three newforms of weight $12$,
level $5$, and trivial character.  The first has rational Fourier
coefficients and $q$-expansion
\[f_{1} = q + 34q^{2} - 792q^{3} - 892q^{4} + 3125q^{5} - 26928q^{6}
-17556q^{7} + \cdots\]
while the other two,
\begin{multline*}
f_{2} = q + (-10+6\sqrt{151})q^{2} + (-110+32\sqrt{151})q^{3} +
(3448-120\sqrt{151})q^{4} \\
- 3125q^{5} +(30092-980\sqrt{151})q^{6} +
(28950+1056\sqrt{151})q^{7} + \cdots
\end{multline*}
and its Galois conjugate, have field of Fourier coefficients
$\Q(\sqrt{151})$.  Note that $\Obs(\bar{f}_{2})$ is simply the set of
conjugates of elements of $\Obs(f_{2})$, so that it
suffices to study $f_{1}$ and $f_{2}$.

Using Lemma~\ref{lemma:red}, one computes that:
\begin{gather*}
\Abs(f_{1}) = \{2,5,31\}; \\
\Abs(f_{2}) = \bigl\{\p_{2},
\p_{5},\bar{\p}_{5},(601,358+\sqrt{151})\bigr\}.
\end{gather*}

We now consider congruences.  By comparing Fourier coefficients, one
sees that $f_{1}$ and $f_{2}$ are congruent modulo primes above $2$ and $5$:
\begin{gather*}
\Cong_{5}(f_{1}) = \{2,5\}; \\
\Cong_{5}(f_{2}) = \bigl\{\p_{2},(5,1+\sqrt{151})\bigr\}.
\end{gather*}
The only possible proper
congruences is with the unique newform
\[\Delta = q -24q^{2} +252q^{3}-1472q^{4}+4830q^{5}-6048q^{6} -16744q^{7}
+ \cdots\]
of weight $12$ and level $1$; one computes that:
\begin{gather*}
\Cong_{<5}(f_{1}) = \{2,29\}; \\
\Cong_{<5}(f_{2}) = \bigl\{\p_{2},(5,4+\sqrt{151}),
(131,46+\sqrt{151})\bigr\}.
\end{gather*}

Both $f_{1}$ and $f_{2}$ are ordinary at $13$, so that by
Theorems~\ref{thm:van} and~\ref{thm:vanc}
we conclude that:
\begin{gather*}
\{29\} \subseteq \Obs(f_{1}) \subseteq \{3,7,11,13,29\}; \\
\bigl\{(131,46+\sqrt{151})\bigr\} \subseteq 
\Obs(f_{2}) \subseteq \bigl\{\p_{3},\bar{\p}_{3},\p_{7},\bar{\p}_{7},
(11),(13),(131,46+\sqrt{151})\bigr\}.
\end{gather*}

\subsection{Weight $\boldsymbol{6}$, level $\boldsymbol{30}$,
trivial character}

There are two newforms of weight $6$, level $30$, and
trivial character, both with rational Fourier coefficients:
\begin{multline*}
f_{1} = q+4q^{2}+9q^{3}-16q^{4}+25q^{5}+36q^{6}+32q^{7}
-192q^{8}\\-162q^{9}+100q^{10}+12q^{11}-144q^{12}-154q^{13} + \cdots
\end{multline*}
\begin{multline*}
f_{2} = q-4q^{2}+9q^{3}-16q^{4}-25q^{5}-36q^{6}+164q^{7}+
192q^{8}\\-162q^{9}+100q^{10}+720q^{11}-144q^{12}+698q^{13} + \cdots
\end{multline*}
Using Lemma~\ref{lemma:red}, one computes that:
\begin{gather*}
\Abs(f_{1}) = \{2,3,5\}; \\
\{2,3\} \subseteq \Abs(f_{2}) \subseteq \{2,3,5\}.
\end{gather*}

The newforms $f_{1}$ and $f_{2}$ have a congruence modulo $12$
(remember that one only checks the Fourier coefficients with
exponent prime to $30$), so that:
\[\Cong_{30}(f_{1}) = \Cong_{30}(f_{2}) = \{2,3\}.\]
There are ten newforms of level dividing $30$ and trivial character
to consider for proper congruences.  The most interesting occur
for $f_{2}$: it has a congruence modulo $19$ with the newform
\begin{multline*}
 q+7q^{2}+9q^{3}+17q^{4}-25q^{5}+63q^{6}+12q^{7}-105q^{8} \\ -162q^{9}
-175q^{10} + 112q^{11} + 153q^{12} - 974q^{13} + \cdots
\end{multline*}
of level $15$, and modulo $31$ with the newform
\begin{multline*}
q-4q^{2}-26q^{3}-16q^{4}-25q^{5}+104q^{6}-22q^{7}+192q^{8}\\
+433q^{9} +100q^{10}-768q^{11}+416q^{12}-46q^{13}+\cdots
\end{multline*}
of level $10$.  In any event, one computes:
\[\Cong_{<30}(f_{1}) = \{2,3,5\}; \quad \Cong_{<30}(f_{2}) = \{2,3,19,31\}. \]

Both $f_{1}$ and $f_{2}$ are ordinary at $7$, so that we conclude that:
\begin{gather*}
\Obs(f_{1}) \subseteq \{7\}; \\
\{19,31\} \subseteq \Obs(f_{2}) \subseteq \{5,7,19,31\}.
\end{gather*}
To give an explicit example of an obstructed set, 
using (\ref{m3}) of Theorem~\ref{thm:vanc} one finds that
\[ \Obs_{\{2,3,5,17,\infty\}}(f_{1}) = \{7\}. \]

\subsection{Weight $\boldsymbol{3}$, level $\boldsymbol{35}$,
character of conductor $\boldsymbol{7}$}

There are four newforms of weight $3$, level $35$, and with quadratic character
\[\omega : (\Z/35\Z)^{\times} \to \{\pm 1\}\]
the Legendre symbol $\left(\frac{\cdot}{7}\right)$.
All four are defined over $\Q(\sqrt{-5})$: two are
\[f_{1} = q - q^{2} - 2\sqrt{-5}q^{3} -3q^{4} -\sqrt{-5}q^{5} + 2\sqrt{-5}q^{6}
+7q^{7}+ \cdots\]
and its Galois conjugate while the other two are
\[f_{2} = q +2q^{2} +\sqrt{-5}q^{3}-\sqrt{-5}q^{5}+2\sqrt{-5}q^{6}
-(2+3\sqrt{-5})q^{7} + \cdots \]
and its Galois conjugate.  As before, it suffices to study $f_{1}$ and $f_{2}$.
Using Lemma~\ref{lemma:red}, one finds that:
\begin{gather*}
\bigl\{\p_{2},\p_{3},\bar{\p}_{3},
\p_{5}\bigr\} \subseteq \Abs(f_{1}) \subseteq
\bigl\{\p_{2},\p_{3},\bar{\p}_{3},\p_{5},\p_{7},\bar{\p}_{7} \bigr\} \\
\bigr\{\p_{3},\bar{\p}_{3}\bigr\} \subseteq
\Abs(f_{2}) \subseteq
\bigr\{\p_{3},\bar{\p}_{3},\p_{5},\p_{7},\bar{\p}_{7}\bigr\}.
\end{gather*}

The newforms $f_{1}$ and $f_{2}$ have a congruence modulo $3$, while
$f_{1}$ and $\bar{f}_{2}$ have no congruences; thus:
\[\Cong_{35}(f_{1}) = \Cong_{35}(f_{2}) = \bigl\{\p_{3},\bar{\p}_{3}\}.\]
Since $\omega$ has conductor $7$,
the only proper congruences we need to check are with the
unique newform
\[q-3q^{2}+5q^{4}-7q^{7}-3q^{8}-9q^{9} -6q^{11} + \cdots\]
of weight $3$, level $7$, and character $\left(\frac{\cdot}{7}\right)$.
One finds that:
\[\Cong_{<35}(f_{1}) = \bigl\{\p_{2}\bigr\};\]
\[\Cong_{<35}(f_{2}) = \bigl\{\p_{5}\bigr\}.\]
Theorem~\ref{thm:van} allows us to conclude that:
\[\Obs(f_{1}) \subseteq \bigl\{\p_{7},\bar{\p}_{7}\bigr\}\]
\[\Obs(f_{2}) \subseteq \bigl\{\p_{2},\p_{5},\p_{7},\bar{\p}_{7}
\bigr\}.\]

\providecommand{\bysame}{\leavevmode\hbox to3em{\hrulefill}\thinspace}

\end{document}